\documentclass[12pt]{article}
\usepackage{latexsym,amssymb,amsmath}
\begin{document}
\baselineskip=20pt

\newcommand{\la}{\langle}
\newcommand{\ra}{\rangle}
\newcommand{\psp}{\vspace{0.4cm}}
\newcommand{\pse}{\vspace{0.2cm}}
\newcommand{\ptl}{\partial}
\newcommand{\dlt}{\delta}
\newcommand{\sgm}{\sigma}
\newcommand{\al}{\alpha}
\newcommand{\be}{\beta}
\newcommand{\G}{\Gamma}
\newcommand{\gm}{\gamma}
\newcommand{\vs}{\varsigma}
\newcommand{\Lmd}{\Lambda}
\newcommand{\lmd}{\lambda}
\newcommand{\td}{\tilde}
\newcommand{\vf}{\varphi}
\newcommand{\yt}{X^{\nu}}
\newcommand{\wt}{\mbox{wt}\:}
\newcommand{\rd}{\mbox{Res}}
\newcommand{\ad}{\mbox{ad}}
\newcommand{\stl}{\stackrel}
\newcommand{\ol}{\overline}
\newcommand{\ul}{\underline}
\newcommand{\es}{\epsilon}
\newcommand{\dmd}{\diamond}
\newcommand{\clt}{\clubsuit}
\newcommand{\vt}{\vartheta}
\newcommand{\ves}{\varepsilon}
\newcommand{\dg}{\dagger}
\newcommand{\tr}{\mbox{Tr}}
\newcommand{\ga}{{\cal G}({\cal A})}
\newcommand{\hga}{\hat{\cal G}({\cal A})}
\newcommand{\Edo}{\mbox{End}\:}
\newcommand{\for}{\mbox{for}}
\newcommand{\kn}{\mbox{ker}}
\newcommand{\Dlt}{\Delta}
\newcommand{\rad}{\mbox{Rad}}
\newcommand{\rta}{\rightarrow}
\newcommand{\mbb}{\mathbb}
\newcommand{\lra}{\Longrightarrow}

\begin{center}{\Large \bf Differential-Operator Representations of
$S_n$}\end{center}
\begin{center}{\Large \bf and Singular Vectors in Verma Modules}  \end{center} \pse

\begin{center}{\large Xiaoping Xu}\end{center}
\begin{center}{Hua Loo-Keng Mathematical Laboratory}\end{center}
\begin{center}{Institute of Mathematics, Academy of Mathematics \& System Sciences}\end{center}
\begin{center}{Chinese Academy of Sciences, Beijing 100190, P.R.China}\footnote{Research supported by CNSF Grant
10871193}\end{center}

\vspace{0.3cm}

\begin{abstract}{ Given a weight of $sl(n,\mbb{C})$, we derive a system of
variable-coefficient second-order linear partial differential
equations that determines the singular vectors in the corresponding
 Verma module, and a differential-operator representation of the symmetric group $S_n$ on
 the related space of truncated power series. We prove that the solution
 space of the system of partial differential equations is exactly
 spanned by $\{\sgm(1)\mid \sgm\in S_n\}$. Moreover, the singular
 vectors of $sl(n,\mbb{C})$ in the Verma module are given by those
 $\sgm(1)$ that are polynomials. The well-known
 results of Verma, Bernstein-Gel'fand-Gel'fand and Jantzen for the case of $sl(n,\mbb{C})$ are
 naturally included in our almost elementary approach of partial differential
 equations.}\end{abstract}

\section{Introduction}

One of the most beautiful things in Lie algebras is the highest
weight representation theory. It was established based on the
induced modules of  a Lie algebra with respect to a Cartan
decomposition  from one-dimensional modules  of the  Borel
subalgebra associated with a linear
 function (weight) on the Cartan subalgebra. These modules are now known as {\it Verma modules}
  [V1]. A {\it singular vector} (or {\it canonical vector}) in a Verma module is
  a weight vector annihilated by positive root vectors. It is well known that the structure of a
  Verma module of a finite-dimensional simple Lie algebra is completely determined by its
  singular vectors (cf. [V1]). In this paper, we find  explicit formulas for singular vectors
   in  Verma modules for the Lie algebra $sl(n,\mbb{C})$ in terms of a differential-operator representation of
   the symmetric group $S_n$ on a certain space of truncated power series.

 The structure of Verma module was first studied by Verma [V1].  Verma reduced the problem of
 determining all submodules of a Verma module of a finite-dimensional  semisimple Lie algebra
  to determining the embeddings of the other Verma modules into the objective module. He
  proved that the multipicity of the embedding is at most one. Bernstein, Gel'fand and
  Gel'fand [BGG] introduced the well-known useful notion of category ${\cal O}$ of
  representations, and found a necessary and sufficient condition for the existence of such a
  embedding in terms of the action of  Weyl group on weights. Sapovolov [S] introduced a certain
   bilinear form on a universal envelopping algebra. Lepowsky [L1-L4] studied analogous induced
    modules with respect to Iwasawa decomposition that is more general than Cartan decomposition,
     and obtained similar results as those in [V1] and [BGG]. These modules are now known as
  {\it generalized Verma modules.}

 Jantzen [J1, J2] introduced his famous ``Jantzen filtrations" on Verma modules and used
Sapovolov form to determine weights of singular vectors in Verma
modules.
 Verma modules of infinite-dimensional Lie algebras were first studied by Kac [Kv1]. Kac and
  Kazhdan [KK] generalized the results of Verma [V1] and Bernstein-Gel'fand-Gel'fand [BGG]  to
   the contragredient Lie algebra corresponding a symmetrizable generalized Cartan matrix.
Deodhar, Gabber and Kac [DGK] generalized the results further to
more general matrix Lie algebras. Rocha-Caridi and Wallach [RW1,
RW2] generalized the results of Verma [V1] and
Bernstein-Gel'fand-Gel'fand [BGG] to a class of graded Lie algebras
possessing a Cartan decomposition and obtained Jantzen's character
formula corresponding to the quotient of two Verma modules. The
resolutions of irreducible highest weight modules over rank-2
Kac-Moody algebras were constructed.

     One of the fundamental and difficult remaining problems in this direction is how to
     determine the singular vectors explicitly. Malikov, Feigin and Fuchs [MFF] introduced a
formal manipulation on products of several general powers of
negative simple root vectors and
 used free Lie algebras to give a rough condition when such product is well defined. It seems
 to us that their condition can not be verified in general and their method can practically be
 applied only to finding very special singular vectors.

     In this paper, we introduce an almost elementary partial differential equation approach
 of determining the singular vectors in any Verma module of $sl(n,\mbb{C})$. First, we identify the
 Verma modules with a space of polynomials, and the action of $sl(n,\mbb{C})$ on the Verma module is
 identified with a differential operator action of $sl(n,\mbb{C})$ on the polynomials. Any singular
 vector in the Verma module becomes a polynomial solution of a system of variable-coefficient
  second-order linear partial differential equations. Thus we have changed a difficult problem
  in a noncommutative space to a problem in  commutative space. However, it is in general impossible to
  solve the system in the space of polynomials. So we extend the action of $sl(n,\mbb{C})$ on the
  polynomial space to a larger space of certain truncated formal power series. On this larger space, the
   negative simple root vectors become differential
operators whose arbitrary complex powers are well defined (so are
their products). In this way,
 we overcome the difficulty of determining whether a product of several general powers of
 negative simple root vectors is well defined in the work [MFF] of Malikov, Feigin and Fuchs.
Next we define a differential-operator representation of
   the symmetric group $S_n$ on the space of truncated power series.
Using commutator relations among root vectors and a certain
substitution-of-variable technique
 that we developed in [X], we prove that  the solution
 space of the system of partial differential equations in the space of truncated power series is exactly
 spanned by $\{\sgm(1)\mid \sgm\in S_n\}$. Moreover, the singular
 vectors of $sl(n,\mbb{C})$ in the Verma module are given by those
 $\sgm(1)$ that are polynomials. In particular, there are exactly
 $n!$ singular vectors up to scalar multiples in the Verma module
 when the weight is dominant integral.

In Section 2, we  derive the system of partial differential
equations and  a differential-operator representation of
   the symmetric group $S_n$ on the space of certain truncated formal power series.
   Moreover, we prove that $\{\sgm(1)\mid \sgm\in S_n\}$ are
the solutions of the system. In Section 3, we  completely solve the
system in the space of power series.

\section{Differential Equations and Representations}

In this section, we first derive a system of variable-coefficient
second-order partial differential equations that determines the
singular vectors in the Verma modules over the special linear Lie
algebra ${\rm sl}(n,\mbb{C})$. Then we construct a differential
operator representation of the symmetric group $S_n$ on the related
space of truncated formal power series and prove that $\{\sgm(1)\mid
\sgm\in S_n\}$ are the solutions of the system in the space.

Denote by $E_{i,j}$ the square matrix with 1 as its $(i,j)$-entry
and 0 as the others.  The special linear Lie algebra
$$sl(n,\mbb{C})=\sum_{1\leq i<j\leq
n}(\mbb{C}E_{i,j}+\mbb{C}E_{j,i})+\sum_{r=1}^{n-1}\mbb{C}(E_{r,r}-E_{r+1,r+1})\eqno(2.1)$$
with the Lie bracket:
$$[A,B]=AB-BA\qquad \for\;\;A,B\in sl(n,\mbb{C}).\eqno(2.2)$$
 Set
$$h_i=E_{i,i}-E_{i+1,i+1},\qquad i=1,2,...,n-1.\eqno(2.3)$$
The subspace
$$H=\sum_{i=1}^{n-1}\mbb{C}h_i\eqno(2.4)$$
forms a Cartan subalgebra of $sl(n,\mbb{C})$. We choose
$$\{E_{i,j}\mid 1\leq i<j\leq n\}\;\;\mbox{as positive root vectors}.\eqno(2.5)$$
In particular, we have
$$\{E_{i,i+1}\mid i=1,2,...,n-1\}\;\;\mbox{as positive simple root vectors}.\eqno(2.6)$$
Accordingly,
$$\{E_{i,j}\mid 1\leq j<i\leq n\}\;\;\mbox{are negative root vectors}\eqno(2.7)$$
and we have
$$\{E_{i+1,i}\mid i=1,2,...,n-1\}\;\;\mbox{as negative simple root vectors}.\eqno(2.8)$$

Denote by $\mbb{N}$ the additive
 semigroup of nonnegative integers.
Let
$$\G=\sum_{1\leq j<i\leq n}\mbb{N}\es_{i,j}\eqno(2.9)$$
be the torsion-free additive semigroup  of  rank $n(n-1)/2$ with
$\es_{i,j}$ as base elements. Let ${\cal G}_-$ be the Lie subalgebra
spanned by (2.7) and let $U({\cal G}_-)$ be its universal enveloping
algebra. For
$$\al =\sum_{1\leq j<i\leq n}\al_{i,j}\es_{i,j}\in \G,\eqno(2.10)$$
we denote
$$E^\al=E_{2,1}^{\al_{2,1}}E_{3,1}^{\al_{3,1}}E_{3,2}^{\al_{3,2}}E_{4,1}^{\al_{4,1}}
\cdots E_{n,1}^{\al_{n,1}}\cdots E_{n,n-1}^{\al_{n,n-1}}\in U({\cal
G}_-).\eqno(2.11)$$
 Then
$$\{E^\al\mid \al\in\G\}\;\;\mbox{forms a basis of}\;\;U({\cal G}_-).\eqno(2.12)$$

 Let $\lmd$ be a weight, which is a linear function
on $H$, such that
$$\lmd(h_i)=\lmd_i\qquad\for\;\;i=1,2,...,n-1.\eqno(2.13)$$
Recall that $sl(n,\mbb{C})$ is generated by
$\{E_{i,i+1},E_{i+1,i}\mid i=1,2,...,n-1\}$ as a Lie algebra. The
Verma  $sl(n,\mbb{C})$-module with the highest-weight vector
$v_\lmd$ of weight $\lmd$ is given by
$$M_\lmd =\mbox{Span}\{E^\al v_\lmd\mid\al\in \G\},\eqno(2.14)$$
with the action determined by
\begin{eqnarray*}E_{i,i+1}(E^\al v_\lmd)&=&(\sum_{j=1}^{i-1}\al_{i+1,j}E^{\al+\es_{i,j}
-\es_{i+1,j}}-\sum_{j=i+2}^n\al_{j,i}E^{\al+\es_{j,i+1}-\es_{j,i}}\\&
&+\al_{i+1,i}(\lmd_i+1
-\sum_{j=i+1}^n\al_{j,i}+\sum_{j=i+2}^n\al_{j,i+1})E^{\al-\es_{i+1,i}})v_\lmd,\hspace{2.3cm}
(2.15)\end{eqnarray*}
$$E_{i+1,i}(E^\al v_\lmd)=(E^{\al+\es_{i+1,i}}+\sum_{j=1}^{i-1}\al_{i,j}E^{\al+\es_{i+1,j}
-\es_{i,j}})v_\lmd\eqno(2.16)$$ for $i=1,...,n-1$. For any
$\al\in\G$, we define the {\it weight} of $E^\al v_\lmd $  by
$$({\rm wt}\:E^\al v_\lmd
)(h_i)=(\lmd_i+\sum_{p=1}^{i-1}(\al_{i,p}-\al_{i+1,p})+\sum_{j=i+2}^n(\al_{j,i+1}-\al_{j,i})
-2\al_{i+1,i})h_i\eqno(2.17)$$
 for $i=1,...,n-1$. Then the Verma module $M_\lmd$ is a space graded by weights.
 A {\it singular vector} is a homogeneous nonzero vector $u$ in $M_\lmd$
such that
$$E_{i,i+1}(u)=0\qquad \for\;\;i=1,...,n-1.\eqno(2.18)$$
Here we have used the fact that all positive root vectors are
generated by simple positive root
 vectors. The Verma module is irreducible if and only if any singular vector is a scalar
 multiple of $v_\lmd$.

Consider the polynomial algebra
$${\cal A}=\mbb{C}[x_{i,j}\mid 1\leq j<i\leq n]\eqno(2.19)$$
in $n(n-1)/2$ variables. Set
$$x^\al=\prod_{1\leq j<i\leq n}x_{i,j}^{\al_{i,j}}\qquad\for\;\;\al\in\G.\eqno(2.20)$$
Then
$$\{x^\al\mid \al\in\G\}\;\;\mbox{forms a basis of}\;\;{\cal A}.\eqno(2.21)$$
Thus we have a linear isomorphism $\tau: M_{\lmd}\rta {\cal A}$
determined by
$$\tau(E^\al v_\lmd)=x^\al\qquad\for\;\;\al\in\G.\eqno(2.22)$$
The algebra ${\cal A}$ becomes $sl(n,\mbb{C})$-module by the action
$$A(f)=\tau(A(\tau^{-1}(f)))\qquad\for\;\;A\in sl(n,\mbb{C}),\;f\in {\cal A}.\eqno(2.23)$$
For convenience, we denote the partial derivatives
$$\ptl_{i,j}=\ptl_{x_{i,j}}\qquad\for\;\;1\leq j<i\leq n.\eqno(2.24)$$
In particular,
\begin{eqnarray*}& &d_i=E_{i,i+1}|_{\cal A}\\ &=&(\lmd_i-\sum_{j=i+1}^nx_{j,i}\ptl_{j,i}
+\sum_{j=i+2}^nx_{j,i+1}\ptl_{j,i+1})\ptl_{i+1,i}+\sum_{j=1}^{i-1}x_{i,j}\ptl_{i+1,j}
-\sum_{j=i+2}^nx_{j,i+1}\ptl_{j,i}\hspace{0.9cm}(2.25)\end{eqnarray*}
for $i=1,2,...,n-1$ by (2.15). \psp

{\bf Proposition 2.1}. {\it A homogeneous vector} $u\in M_\lmd$ {\it
is a singular vector if and only if}
$$d_i(\tau(u))=0\qquad\for\;\;i=1,2,...,n-1.\eqno(2.26)$$
\pse

The system of partial differential equations
\begin{eqnarray*}\hspace{2.5cm}& &(\lmd_i-\sum_{j=i+1}^nx_{j,i}\ptl_{j,i}+\sum_{j=i+2}^n
x_{j,i+1}\ptl_{j,i+1})\ptl_{i+1,i}(z)\\ &
&+\sum_{j=1}^{i-1}x_{i,j}\ptl_{i+1,j}(z)-
\sum_{j=i+2}^nx_{j,i+1}\ptl_{j,i}(z)=0\hspace{4.3cm}(2.27)\end{eqnarray*}
for $i=1,2,...,n-1$ and unknown function $z$ in $\{x_{i,j}\mid 1\leq
j<i\leq n\}$,  is called
 the {\it system of partial differential equations for the singular vectors of $sl(n,\mbb{C})$}.

Next we want to construct a differential operator representation of
the symmetric group $S_n$ on the related space of truncated series
and prove that $\{\sgm(1)\mid \sgm\in S_n\}$ are the solutions of
the system.
 First, we have
$$\eta_i=E_{i+1,i}|_{\cal A}=x_{i+1,i}+\sum_{j=1}^{i-1}x_{i+1,j}\ptl_{i,j}\eqno(2.28)$$
for $i=1,2,...,n-1$ by (2.16). Now we view $\{d_i,\eta_i\mid
i=1,2,...,n-1\}$ purely as differential operators acting on
functions of $\{x_{i,j}\mid 1\leq j<i\leq n\}$. In this way,
 we get a Lie algebra action on functions of $\{x_{i,j}\mid 1\leq j<i\leq n\}$ through
 $E_{i,i+1}=d_i$ and $E_{i+1,i}=\eta_i$ because
$sl(n,\mbb{C})$ is generated by $\{E_{i,i+1},E_{i+1,i}\mid
i=1,2,...,n-1\}$ as a Lie algebra.  Note
 that
$$h_i(E^\al v_\lmd)=(\lmd_i+\sum_{p=1}^{i-1}(\al_{i,p}-\al_{i+1,p})+\sum_{j=i+2}^n(\al_{j,i+1}
-\al_{j,i})-2\al_{i+1,i})E^\al v_\lmd\eqno(2.29)$$ for
$i=1,2,...,n-1$ and $\al\in\G$. Accordingly, we set
$$\zeta_i=h_i|_{\cal A}=\lmd_i+\sum_{p=1}^{i-1}(x_{i,p}\ptl_{i,p}-x_{i+1,p}\ptl_{i+1,p})
+\sum_{j=i+2}^n(x_{j,i+1}\ptl_{j,i+1}-x_{j,i}\ptl_{j,i})-2x_{i+1,i}\ptl_{i+1,i}\eqno(2.30)$$
for $i=1,2,...,n-1$. The elements $h_i$ act on functions of
$\{x_{i,j}\mid 1\leq j<i\leq n\}$ through $\zeta_i$.  A function $f$
of $\{x_{i,j}\mid 1\leq j<i\leq n\}$ is called {\it weighted} if
there exist constants $\mu_1,\;\mu_2,\;...,\;\mu_{n-1}$ such that
$$\zeta_i(f)=\mu_if\qquad\for\;\;i=1,2,...,n-1.\eqno(2.31)$$
Since $d_i$ maps weighted functions to weighted functions, the
system (2.27) is a weighted system. Any nonzero weighted solution of
the system (2.27) is a singular vector of $sl(n,\mbb{C})$. In
particular, any nonzero  weighted  polynomial solution $f$ of the
system (2.27) gives a singular vector $\tau^{-1}(f)$ in the Verma
module $M_{\lmd}$.

Let
$${\cal A}_0=\mbb{C}[x_{i,j}\mid 1\leq j\leq i-2\leq n-2]\eqno(2.32)$$
be the polynomial algebra in $\{x_{i,j}\mid 1\leq j\leq i-2\leq
n-2\}$. We denote
$$x^{\vec a}=\prod_{i=1}^{n-1}x_{i+1,i}^{a_i}\qquad\for\;\;\vec a=(a_1,a_2,...,
a_{n-1})\in\mbb{C}^{\:n-1}.\eqno(2.33)$$ Let
$${\cal A}_1=\{\sum_{\vec j\in\mbb{N}^{\;n-1}}\sum_{i=1}^pf_{\vec a\:^i-\vec j}
x^{\vec a\:^i-\vec j}\mid 1\leq p\in\mbb{N},\;\vec a\:^i\in
\mbb{C}^{\:n-1}, \;f_{\vec a\:^i-\vec j}\in{\cal A}_0\}\eqno(2.34)$$
be the space of truncated-up formal power series in
$\{x_{2,1},x_{3,2},...,x_{n,n-1}\}$ over ${\cal A}_0$. Then ${\cal
A}$ is a subspace of ${\cal A}_1$. Since ${\cal A}_1$ is invariant
under the action of $\{E_{i,i+1}|_{{\cal A}_1}=d_i,
E_{i+1,i}|_{{\cal A}_1}=\eta_i\mid i=1,2,...,n-1\}$, ${\cal A}_1$
becomes an $sl(n,\mbb{C})$-module.

For $a\in \mbb{C}$ and $p\in\mbb{N}$, we denote
$$\la a\ra_p=a(a-1)(a-2)\cdots (a-p+1).\eqno(2.35)$$
 Moreover, by (2.28), we define
$$\eta^{a}_i=(x_{i+1,i}+\sum_{j=1}^{i-1}x_{i+1,j}\ptl_{i,j})^{a}=\sum_{p=0}^{\infty}
\frac{\la
a\ra_p}{p!}x_{i+1,i}^{a-p}(\sum_{j=1}^{i-1}x_{i+1,j}\ptl_{i,j})^p\eqno(2.36)$$
as differential operators on ${\cal A}_1$, for $i=1,2,...,n-1$ and
$a\in\mbb{C}$. If $a\not\in\mbb{N}$, then the above summation is
infinite  and the positions of $x_{i+1,i}$ and
$(\sum_{j=1}^{i-1}x_{i+1,j}\ptl_{i,j})$ are not symmetric. Since
$x_{i+1,i}$ and $(\sum_{j=1}^{i-1}x_{i+1,j}\ptl_{i,j})$ commute, we
have
$$\eta_i^{a_1}\eta_i^{a_2}=\eta^{a_1+a_2}_i\qquad\for\;\;a_1,a_2\in\mbb{C}.
\eqno(2.37)$$ In particular, the inverse of the differential
operator $\eta_i^{a}$ is exactly  $\eta_i^{-a}$.

Given two differential operators $d$ and $\bar{d}$, we define the
commutator
$$[d,\bar{d}]=d\bar{d}-\bar{d}d.\eqno(2.38)$$
 For any element $f\in {\cal A}_1$ and $r\in\mbb{C}$, we have
$$[\ptl_{i+1,i},x_{i+1,i}^r](f)=\ptl_{i+1,i}(x_{i+1,i}^rf)-x_{i+1,i}^r\ptl_{i+1,i}(f)=
rx_{i+1,i}^{r-1}f,\eqno(2.39)$$ that is,
$$[\ptl_{i+1,i},x_{i+1,i}^r]=rx_{i+1,i}^{r-1}\qquad\mbox{as operators}.\eqno(2.40)$$
Note that if $(r,s)\not\in\{(i+1,j)\mid j=1,...,i\}$ and
$(p,q)\not\in\{(i,j)\mid
 j=1,...,i-1\}$,
 then
$$[\ptl_{r,s},\eta^{a}_i]=[x_{p,q},\eta^{a}_i]=0\eqno(2.41)$$
directly by (2.36). Now
\begin{eqnarray*}\hspace{1cm}[\ptl_{i+1,i},\eta^{a}_i] &=&\sum_{p=0}^{\infty}
\frac{\la
a\ra_p}{p!}[\ptl_{i+1,i},x_{i+1,i}^{a-p}(\sum_{j=1}^{i-1}x_{i+1,j}\ptl_{i,j})^p]
\\ &=&\sum_{p=0}^{\infty}\frac{\la a\ra_p}{p!}(a-p)x_{i+1,i}^{a-p-1}
(\sum_{j=1}^{i-1}x_{i+1,j}\ptl_{i,j})^p\\
&=&\sum_{p=0}^{\infty}\frac{a\la a-1\ra_p}{p!}x_{i+1,i}^{a-p-1}
(\sum_{j=1}^{i-1}x_{i+1,j}\ptl_{i,j})^p=a
\eta_i^{a-1}\hspace{3.2cm}(2.42)\end{eqnarray*} by (2.40). Moreover,
for $j=1,2,...,i-1$,
\begin{eqnarray*}\hspace{1cm}[\ptl_{i+1,j},\eta^{a}_i]&=&\sum_{p=0}^{\infty}
\frac{\la
a\ra_p}{p!}[\ptl_{i+1,j},x_{i+1,i}^{a-p}(\sum_{s=1}^{i-1}x_{i+1,s}\ptl_{i,s})^p]
\\ &=& \sum_{p=0}^{\infty}
\frac{\la
a\ra_p}{p!}px_{i+1,i}^{a-p}(\sum_{s=1}^{i-1}x_{i+1,s}\ptl_{i,s})^{p-1}\ptl_{i,j}
\\ &=& \sum_{p=0}^{\infty}
\frac{a\la
a-1\ra_{p-1}}{(p-1)!}x_{i+1,i}^{a-p}(\sum_{s=1}^{i-1}x_{i+1,s}\ptl_{i,s})^{p-1}
\ptl_{i,j}=a
\eta_i^{a-1}\ptl_{i,j}\hspace{1.5cm}(2.43)\end{eqnarray*} and
similarly,
$$[x_{i,j},\eta^{a}_i]=-a \eta^{a-1}_ix_{i+1,j}\qquad\for\;\;j=1,2,...,i-1.\eqno(2.44)$$
\vspace{0.1cm}

{\bf Lemma 2.2}. {\it For} $i,l\in\{1,2,...,n-1\}$ {\it and}
$a\in\mbb{C}$, {\it we have}:
$$[d_l,\eta_i^{a}]=a\dlt_{i,l} \eta_i^{a-1}(1-a+\zeta_i).\eqno(2.45)$$

{\it Proof}. Note that
$$[E_{l,l+1},E_{i+1,i}^m]=m\dlt_{i,l}E_{i+1,i}^{m-1}(1-m+h_i)\qquad\for\;\;m\in\mbb{N}
\eqno(2.46)$$ (cf. (2.3)). So (2.45) holds for any $a\in\mbb{N}$ by
(2.25), (2.28) and (2.30). Since (2.45)
 is completely determined by  (2.41)-(2.44), which are independent of whether $a$ is a
 nonnegative integer, it must hold for any $a\in\mbb{C}.\qquad\Box$
\psp

Denote the  Cartan matrix of $sl(n,\mbb{C})$ by
$$\left[\begin{array}{cccc}a_{1,1}&a_{1,2}&\cdots& a_{1,{n-1}}\\ a_{2,1}&a_{2,2}&\cdots
&a_{2,{n-1}}\\ \vdots&\vdots& &\vdots\\ a_{n-1,1}&a_{n-1,2}&\cdots&
a_{n-1,{n-1}}
\end{array}\right]=\left[\begin{array}{rrrr}2&-1&&\\ -1&2&\ddots &\\ &\ddots& \ddots&-1
\\ &&-1&2\end{array}\right].\eqno(2.47)$$
\vspace{0.1cm}

{\bf Lemma 2.3}. {\it For} $i,l\in\{1,2,...,n-1\}$ {\it and}
$a\in\mbb{C}$,
$$[\zeta_l,\eta_i^{a}]=-a a_{l,i}\eta_i^{a}.\eqno(2.48)$$

{\it Proof}. Observe that
$$[h_l,E_{i+1,i}^m]=ma_{l,i}E_{i+1,i}^m\qquad\for\;\;m\in\mbb{N}\eqno(2.49)$$
(cf. (2.3)). Hence (2.48) holds for any $a\in\mbb{N}$ by  (2.28) and
(2.30). Again  (2.48) is completely determined by (2.41)-(2.44),
which are independent of whether
 $a$ is a nonnegative integer. Thus (2.48) must hold for any $a\in\mbb{C}.\qquad\Box$
\psp

In order to construct a differential-operator representation of the
symmetric group $S_n$ on ${\cal A}_1$, we need the following
result.\psp

{\bf Lemma 2.4}. {\it For any $a_1,a_2\in\mbb{C}$ and $1\leq i<n-1$,
we have}
$$\eta_i^{a_1}\eta_{i+1}^{a_1+a_2}\eta_i^{a_2}
=\eta_{i+1}^{a_2}\eta_i^{a_1+a_2}\eta_{i+1}^{a_1}.\eqno(2.50)$$

{\it Proof}. Note that for $a\in\mbb{C}$, we have
$$[\sum_{p=1}^ix_{i+2,p}\ptl_{i+1,p},x_{i+1,i}^{a}]=a
x_{i+1,i}^{a-1}x_{i+2,i}\eqno(2.51)$$ by (2.40). Moreover,
$$[\sum_{p=1}^ix_{i+2,p}\ptl_{i+1,p},\sum_{j=1}^{i-1}x_{i+1,j}\ptl_{i,j}]=\sum_{j=1}^{i-1}
x_{i+2,j}\ptl_{i,j}.\eqno(2.52)$$ Hence
\begin{eqnarray*} & &\eta_i^{a_1}\eta_{i+1}^{a_1+a_2}\\
&=&\sum_{p,q=0}^{\infty}\frac{\la a_1\ra_p\la
a_1+a_2\ra_q}{p!q!}x_{i+2,i+1}^{a_1+a_2-q}
x_{i+1,i}^{a_1-p}(\sum_{j_1=1}^{i-1}x_{i+1,j_1}\ptl_{i,j_1})^p(\sum_{j_2=1}^ix_{i+2,j_2}
\ptl_{i+1,j_2})^q\\
&=&\sum_{p,q,r,s=0}^{\infty}\frac{(-1)^{r+s}\la a_1\ra_{p+r}\la
a_1+a_2\ra_q\la p\ra_s\la
q\ra_{r+s}}{r!s!p!q!}x_{i+2,i+1}^{a_1+a_2-q}(\sum_{j_2=1}^ix_{i+2,j_2}
\ptl_{i+1,j_2})^{q-r-s}x_{i+1,i}^{a_1-p-r}\\ &
&\times(\sum_{j_1=1}^{i-1}x_{i+1,j_1}\ptl_{i,j_1})^{p-s}x_{i+2,i}^r(\sum_{j=1}^{i-1}
x_{i+2,j}\ptl_{i,j})^s\\
&=&\sum_{p,q,r,s=0}^{\infty}\frac{(-1)^{r+s}\la a_1\ra_{p+r}\la
a_1+a_2\ra_q}{r!s!(p-s)!(q-r-s)!}x_{i+2,i+1}^{a_1+a_2-q}(\sum_{j_2=1}^ix_{i+2,j_2}
\ptl_{i+1,j_2})^{q-r-s}x_{i+1,i}^{a_1-p-r}\\ &
&\times(\sum_{j_1=1}^{i-1}x_{i+1,j_1}\ptl_{i,j_1})^{p-s}x_{i+2,i}^r(\sum_{j=1}^{i-1}
x_{i+2,j}\ptl_{i,j})^s\hspace{5cm}\end{eqnarray*}
\begin{eqnarray*}
&=&\sum_{q,k,s=0}^{\infty}\;\sum_{p=0}^{\infty} \frac{(-1)^k\la
a_1\ra_k\la a_1-k\ra_{p-s}\la a_1+a_2\ra_q}{(k-s)!s!(p-s)!(q-k)!}
x_{i+2,i+1}^{a_1+a_2-q}(\sum_{j_2=1}^ix_{i+2,j_2}
\ptl_{i+1,j_2})^{q-k}\\ & &\times
x_{i+1,i}^{a_1-k-(p-s)}(\sum_{j_1=1}^{i-1}x_{i+1,j_1}\ptl_{i,j_1})^{p-s}x_{i+2,i}^{k-s}
(\sum_{j=1}^{i-1} x_{i+2,j}\ptl_{i,j})^s\\
&=&\sum_{q,k,s=0}^{\infty} \frac{(-1)^k\la a_1\ra_k\la
a_1+a_2\ra_q}{(k-s)!s!(q-k)!}
x_{i+2,i+1}^{a_1+a_2-q}(\sum_{j_2=1}^ix_{i+2,j_2}
\ptl_{i+1,j_2})^{q-k}\\ & &\times \eta_i^{a_1-k}x_{i+2,i}^{k-s}
(\sum_{j=1}^{i-1}
x_{i+2,j}\ptl_{i,j})^s\\&=&\sum_{q,k=0}^{\infty}\;\sum_{s=0}^{\infty}
\frac{(-1)^k\la a_1\ra_k\la a_1+a_2\ra_q}{(k-s)!s!(q-k)!}
x_{i+2,i+1}^{a_1+a_2-q}(\sum_{j_2=1}^ix_{i+2,j_2}
\ptl_{i+1,j_2})^{q-k}\\ & &\times \eta_i^{a_1-k}x_{i+2,i}^{k-s}
(\sum_{j=1}^{i-1} x_{i+2,j}\ptl_{i,j})^s\\
&=&\sum_{q,k=0}^{\infty}\frac{(-1)^k\la a_1\ra_k\la
a_1+a_2\ra_q}{k!(q-k)!}
x_{i+2,i+1}^{a_1+a_2-q}(\sum_{j_2=1}^ix_{i+2,j_2}
\ptl_{i+1,j_2})^{q-k}\\ & &\times \eta_i^{a_1-k}(x_{i+2,i}+
\sum_{j=1}^{i-1} x_{i+2,j}\ptl_{i,j})^k\\
&=&\sum_{k=0}^{\infty}\;\sum_{q=0}^{\infty}\frac{(-1)^k\la
a_1\ra_k\la a_1+a_2\ra_k \la a_1+a_2-k\ra_{q-k} }{k!(q-k)!}
x_{i+2,i+1}^{a_1+a_2-k-(q-k)}(\sum_{j_2=1}^ix_{i+2,j_2}
\ptl_{i+1,j_2})^{q-k}\\ & &\times \eta_i^{a_1-k}(x_{i+2,i}+
\sum_{j=1}^{i-1}
x_{i+2,j}\ptl_{i,j})^k\\&=&\sum_{k=0}^{\infty}\frac{(-1)^k\la
a_1\ra_k\la a_1+a_2\ra_k
}{k!}\eta_{i+1}^{a_1+a_2-k}\eta_i^{a_1-k}(x_{i+2,i}+
\sum_{j=1}^{i-1}
x_{i+2,j}\ptl_{i,j})^k\hspace{2.7cm}(2.53)\end{eqnarray*} by (2.36),
(2.51) and (2.52). Similarly, we have
$$\eta_i^{a_1+a_2}\eta_{i+1}^{a_1}=\sum_{k=0}^{\infty}\frac{(-1)^k\la a_1\ra_k\la a_1
+a_2\ra_k }{k!}\eta_{i+1}^{a_1-k}\eta_i^{a_1+a_2-k}(x_{i+2,i}+
\sum_{j=1}^{i-1} x_{i+2,j}\ptl_{i,j})^k.\eqno(2.54)$$ Thus
\begin{eqnarray*}\eta_i^{a_1}\eta_{i+1}^{a_1+a_2}\eta_i^{a_2}&=&\sum_{k=0}^{\infty}\frac{(-1)^k\la a_1\ra_k\la a_1+a_2\ra_k
}{k!}\eta_{i+1}^{a_1+a_2-k}\eta_i^{a_1+a_2-k}(x_{i+2,i}+
\sum_{j=1}^{i-1} x_{i+2,j}\ptl_{i,j})^k\\
&=&\eta_{i+1}^{a_2}
\eta_i^{a_1+a_2}\eta_{i+1}^{a_1}.\qquad\Box\hspace{7.9cm}(2.55)\end{eqnarray*}
\psp

It is well known that the symmetric group $S_n$ is a group generated
by $\{\sgm_1,...,\sgm_{n-1}\}$ with the defining relations:
$$\sgm_i\sgm_{i+1}\sgm_i=\sgm_{i+1}\sgm_i\sgm_{i+1},\qquad\sgm_r\sgm_s=\sgm_s\sgm_r,\qquad\sgm_r^2=1\eqno(2.56)$$
for $i=1,2,...,n-2$ and $r,s=1,2,...,n-1$ such that $|r-s|\geq 2$.
According to (2.31) and (2.34), any element $f\in {\cal A}_1$ can be
written as $f=\sum_{j\in\mbb{Z}}f_j$ such that
$$\zeta_i(f_j)=\mu_{(j)}(h_i)f_j,\qquad \mu_{(j)}\in
H^\ast.\eqno(2.57)$$ We define an action of
$\{\sgm_1,...,\sgm_{n-1}\}$ on ${\cal A}_1$ by
$$\sgm_i(f)=\sum_{j\in\mbb{Z}}\eta^{\mu_{(j)}(h_i)+1}_i(f_j),\qquad
i=1,2,...,n-1.\eqno(2.58)$$ \pse

{\bf Theorem 2.5}. {\it Expression (2.58) gives a representation of
the symmetric group $S_n$. Moreover, $\{\sgm(1)\mid \sgm\in S_n\}$
are weighted solutions of the system (2.27) of partial differential
equations.}

{\it Proof}. Let $f\in{\cal A}_1$ be a weight function with weight
$\mu$, that is, $\zeta_i(f)=\mu(h_i)f$ for $i=1,2,...,n-1.$ Then
$$(\sgm_i|_{{\cal A}_1})^2(f)=\sgm_i(\eta_i^{\mu(h_i)+1}(f))=\eta_i^{-\mu(h_i)-1}(\eta_i^{\mu(h_i)+1}(f))=f\eqno(2.59)$$
by (2.37), (2.47), (2.48) and (2.58). Thus
$$(\sgm_i|_{{\cal A}_1})^2=\mbox{Id}_{{\cal A}_1}\qquad\for\;\;i=1,2,...,n-1.\eqno(2.60)$$
Note \begin{eqnarray*}& &[(\sgm_i|_{{\cal A}_1})(\sgm_{i+1}|_{{\cal
A}_1})(\sgm_i|_{{\cal
A}_1})](f)\\&=&\sgm_i[\sgm_{i+1}(\eta_i^{\mu(h_i)+1}(f))]=\sgm_i[\eta_{i+1}^{\mu(h_{i+1})+\mu(h_i)+2}
(\eta_i^{\mu(h_i)+1}(f))]\\
&=&\eta_i^{\mu(h_{i+1})+1}[\eta_{i+1}^{\mu(h_{i+1})+\mu(h_i)+2}
(\eta_i^{\mu(h_i)+1}(f))]=(\eta_i^{\mu(h_{i+1})+1}\eta_{i+1}^{\mu(h_{i+1})+\mu(h_i)+2}
\eta_i^{\mu(h_i)+1})(f)\\
&=&(\eta_{i+1}^{\mu(h_i)+1}\eta_i^{\mu(h_{i+1})+\mu(h_i)+2}
\eta_{i+1}^{\mu(h_{i+1})+1})(f)=\eta_{i+1}^{\mu(h_i)+1}[\eta_i^{\mu(h_{i+1})+\mu(h_i)+2}
(\eta_{i+1}^{\mu(h_{i+1})+1}(f))]\\ &=&[(\sgm_{i+1}|_{{\cal
A}_1})(\sgm_i|_{{\cal A}_1})(\sgm_{i+1}|_{{\cal
A}_1})](f)\hspace{8.7cm}(2.61)\end{eqnarray*} by (2.47), (2.48),
(2.50) and (2.58). Hence
$$(\sgm_i|_{{\cal A}_1})(\sgm_{i+1}|_{{\cal
A}_1})(\sgm_i|_{{\cal A}_1})=(\sgm_{i+1}|_{{\cal
A}_1})(\sgm_i|_{{\cal A}_1})(\sgm_{i+1}|_{{\cal A}_1}).\eqno(2.62)$$
For $|r-s|\geq 2$, we have
$$\eta_r^a\eta_s^b=(x_{r+1,r}+\sum_{j=1}^{r-1}x_{r+1,j}\ptl_{r,j})^{a}
(x_{s+1,s}+\sum_{j=1}^{s-1}x_{s+1,j}\ptl_{s,j})^b=\eta_s^b\eta_r^a\eqno(2.63)$$
for $a,b\in\mbb{C}$ by (2.36). So
$$[(\sgm_r|_{{\cal A}_1})(\sgm_s|_{{\cal
A}_1})](f)=(\eta_r^{\mu(h_r)+1}\eta_s^{\mu(h_s)+1})(f)=[(\sgm_s|_{{\cal
A}_1})(\sgm_r|_{{\cal A}_1})](f),\eqno(2.64)$$ which implies
$$(\sgm_r|_{{\cal A}_1})(\sgm_s|_{{\cal
A}_1})=(\sgm_s|_{{\cal A}_1})(\sgm_r|_{{\cal A}_1}).\eqno(2.65)$$
According to (2.56), this proves that (2.58) defines a
representation of $S_n$ on ${\cal A}_1$.

Now we assume that $f\in{\cal A}_1$ is a weighted solution of (2.27)
with weight $\mu$, that is, $d_i(f)=0$ for $i=1,2,...,n-1$. Given
$r\in\{1,2,...,n-1\}$,
\begin{eqnarray*}\hspace{2cm}& &d_i(\sgm_r(f))=d_i\eta_r^{\mu(h_r)+1}(f)=[d_i,\eta_r^{\mu(h_r)+1}](f)
\\ &=&(\mu(h_r)+1)\dlt_{i,r}
\eta_r^{\mu(h_r)}(-\mu(h_r)+\zeta_r)(f)=0\hspace{4.6cm}(2.66)\end{eqnarray*}
by (2.45). So $\sgm_r(f)$ is also a weighted solution of (2.27).
Recall that $S_n$ is generated by $\sgm_1,\sgm_2,...,\sgm_{n-1}$.
For any $\sgm\in S_n$, we write
$\sgm=\sgm_{i_1}\sgm_{i_2}\cdots\sgm_{i_r}$ and
$$\sgm(1)=\sgm_{i_1}(\sgm_{i_2}(\cdots(\sgm_{i_r}(1))\cdots))\eqno(2.67)$$
is a weighted solution of (2.27) by (2.66) and induction on
$r.\qquad\Box$

\section{Completeness}

In this section, we want to prove the following theorem of
completeness: \psp

{\bf Theorem 3.1}. {\it The solution space of the system (2.27) in
${\cal A}_1$ is spanned by $\{\sgm(1)\mid \sgm\in S_n\}$. Moreover,
$\{\sgm(1)\mid \sgm\in S_n\}$ are all the weighted solutions of the
system (2.27) in ${\cal A}_1$ up to scalar multiples. In particular,
there are exactly
 $n!$ singular vectors up to scalar multiples in the Verma module
 $M_\lmd$  when the weight $\lmd$ is dominant integral. }

{\it Proof}. For any
$$\vec a=(a_1,a_2,...,a_{n-1})\in \mbb{C}^{\:n-1},\eqno(3.1)$$
we define
\begin{eqnarray*}\hspace{1.2cm}\phi_{\vec a}&=&\eta_2^{a_2}\eta_1^{a_2-\lmd_{n-1}-1}
\cdots
\eta_i^{a_i}\eta_{i-1}^{a_i-\lmd_{n-1}-1}\cdots\eta_1^{a_i-\lmd_{n-1}-\cdots
-\lmd_{n-(i-1)}-(i-1)}\\ & &\cdots
\eta_{n-1}^{a_{n-1}}\eta_{n-2}^{a_{n-1}-\lmd_{n-1}-1}
\cdots\eta_1^{a_{n-1}+2-\lmd_2-\cdots-\lmd_{n-1}-n}(1).\hspace{3.8cm}(3.2)\end{eqnarray*}
Then
$$\phi_{\vec a}\in {\cal A}_1\eqno(3.3)$$
(cf. (2.32)-(2.34)) and is a solution of the system
$$d_2(z)=d_3(z)=\cdots=d_{n-1}(z)=0\eqno(3.4)$$
by Lemmas 2.2, 2.3 and (2.66). \psp

{\bf Claim}. {\it An element} $z$ {\it  in} ${\cal A}_1$ {\it is a
solution of the system
 (3.4) if and only if it can be written as}
$$z=\sum_{\vec j\in\mbb{N}^{\;n-1}}\sum_{i=1}^pc_{\vec a\:^i-\vec j}x_{2,1}^{a^i_1-j_1}
\phi_{\vec a\:^i-\vec j}\qquad\mbox{\it with}\;\;c_{\vec a\:^i-\vec
j}\in \mbb{C} \eqno(3.5)$$ {\it for some} $\vec a\:^1,...,\vec
a\:^p\in \mbb{C}^{\:n-1}$. \psp

Recall $\eta_1=x_{2,1}$. By Lemma 2.2, the sufficiency holds. Now we
want to prove the necessity. Recall that
$$E_{i+1,i}=d_i,\;\;E_{i,i+1}=\eta_i\qquad\mbox{as operators on}\;\;{\cal A}_1\eqno(3.6)$$
(cf. (2.25) and (2.28)) for $i=1,2,...,n-1$. Note
$$d_{n-1}=(\lmd_{n-1}-x_{n,n-1}\ptl_{n,n-1})\ptl_{n,n-1}+\sum_{i=1}^{n-2}x_{n-1,i}\ptl_{n,i}.
\eqno(3.7)$$ Moreover, (2.11) tells us that
\begin{eqnarray*}\hspace{1.5cm}& &d_{n-2,n}=E_{n-2,n}|_{{\cal A}_1}\\ &=&
(\lmd_{n-1}+\lmd_{n-2}-x_{n,n-2}\ptl_{n,n-2}-x_{n,n-1}\ptl_{n,n-1})\ptl_{n,n-2}\\
& &-d_{n-1}\ptl_{n-1,n-2}+\sum_{i=1}^{n-3}x_{n-2,i}\ptl_{n,i}\\
&=&(\lmd_{n-1}+\lmd_{n-2}
+1-x_{n,n-2}\ptl_{n,n-2}-x_{n,n-1}\ptl_{n,n-1})\ptl_{n,n-2}\\ &
&-\ptl_{n-1,n-2}d_{n-1}
+\sum_{i=1}^{n-3}x_{n-2,i}\ptl_{n,i}.\hspace{7.1cm}(3.8)\end{eqnarray*}
Set
$$\bar{\lmd}_i=n-i-1+\sum_{p=i}^{n-1}\lmd_p\qquad\for\;\;i=2,3,...,n-2.\eqno(3.9)$$
Furthermore,
\begin{eqnarray*}& &d_{n-3,n}=E_{n-3,n}|_{{\cal A}_1}\\ &=& (\bar{\lmd}_{n-3}-2-x_{n,n-3}
\ptl_{n,n-3}-x_{n,n-2}\ptl_{n,n-2}-x_{n,n-1}\ptl_{n,n-1})\ptl_{n,n-3}\\
& &-d_{n-2,n}
\ptl_{n-2,n-3}-d_{n-1}\ptl_{n-1,n-3}+\sum_{i=1}^{n-4}x_{n-3,i}\ptl_{n,i}
\\ &=&(\bar{\lmd}_{n-2}-2-\sum_{p=1}^3
x_{n,n-p}\ptl_{n,n-p}))\ptl_{n,n-3}
-\ptl_{n-2,n-3}d_{n-2,n}-\ptl_{n-1,n-3}d_{n-1}
\\&&-[d_{n-2,n},\ptl_{n-2,n-3}]-[d_{n-1},\ptl_{n-1,n-3}]+\sum_{i=1}^{n-4}x_{n-3,i}\ptl_{n,i}
\\ &=&
(\bar{\lmd}_{n-3}-\sum_{p=1}^3x_{n,n-p}\ptl_{n,n-p})\ptl_{n,n-3}
+\sum_{i=1}^{n-4}x_{n-3,i}\ptl_{n,i}\\ &
&-\ptl_{n-2,n-3}d_{n-2,n}-\ptl_{n-1,n-3}d_{n-1}.\hspace{8.2cm}(3.10)\end{eqnarray*}
By induction, we can prove that
$$d_{i,n}=E_{i,n}|_{{\cal A}_1}=(\bar{\lmd}_i-\sum_{p=i}^{n-1}x_{n,p}\ptl_{n,p})\ptl_{n,i}
+\sum_{q=1}^{i-1}x_{i,q}\ptl_{n,q}-\sum_{j=i+1}^{n-1}\ptl_{j,i}d_{j,n}\eqno(3.11)$$
for $i=2,3,...,n-2$, where we take
$$d_{n-1,n}=d_{n-1}.\eqno(3.12)$$

Suppose that
$$z=fx_{2,1}^{a_1}\phi_{\vec a}\eqno(3.13)$$
is a solution of the system (3.4) for some $\vec a\in
\mbb{C}^{\:n-1}$ and $f\in{\cal A}_0$
 (cf. (2.32) and (3.2)). We want to prove that  $f$ is a constant.  Denote
$$\es_i=(0,...,0,\stl{i}{1},0,...,0)\in\mbb{C}^{\:n-1}\eqno(3.14)$$
and
$$\iota_{i,j}=a_{i+j}-\sum_{p=1}^j(\lmd_{n-p}+1)\qquad\for\;\;2\leq i\leq n-1,\;0\leq j\leq
 n-i-1.\eqno(3.15)$$
We define
$$U_i=\{\sum_{\vec j\in\mbb{N}^{\:n-1}}g_{\vec j}x_{2,1}^{a_1-j_1}\phi_{\vec a-\es_i
-\vec j}\mid g_{\vec j}\in {\cal A}_0\}\eqno(3.16)$$ for
$i=2,...,n-1$, and
$$U=\sum_{i=2}^{n-1}U_i.\eqno(3.17)$$
Moreover, for fixed $i\geq 2$,
$$\{\eta_i^{\iota_{i,j}}\mid 0\leq j\leq n-i-1\}\eqno(3.18)$$
are all the factors in the righthand side of (3.2) that contain
$x_{i+1,p}$ or $\ptl_{i,q}$ with $p=1,...,i$ and $q=1,2,....,i-1$ by
(2.36). Besides,
$$[\ptl_{i+1,i},\eta_i^{\iota_{i,j}}]=\iota_{i,j}\eta_i^{\iota_{i,j}-1},\eqno(3.19)$$
$$[\ptl_{i+1,p},\eta_i^{\iota_{i,j}}]=\iota_{i,j}\eta_i^{\iota_{i,j}-1}\ptl_{i,p}
\qquad\for\;\;p=1,...,i-1,\eqno(3.20)$$
$$[\sum_{p=r}^ix_{i+1,p}\ptl_{i+1,p},\eta_i^{\iota_{i,j}}]=\iota_{i,j}(\eta_i^{\iota_{i,j}}
-\sum_{p=1}^{r-1}x_{i+1,p}\eta_i^{\iota_{i,j}-1}\ptl_{i,p}),\eqno(3.21)$$
$$[\sum_{p=q}^{i-1}x_{i,p}\ptl_{i,p},\eta_i^{\iota_{i,j}}]=-\iota_{i,j}\sum_{p=q}^{i-1}
x_{i+1,p}\eta_i^{\iota_{i,j}-1}\ptl_{i,p},\eqno(3.22)$$ where $2\leq
r\leq i$ and $2\leq q\leq i-1$.

By (3.18)-(3.22), we have
$$\ptl_{i+1,r}(\phi_{\vec a})\in  U_i\eqno(3.23)$$
and
$$(\sum_{p=r}^ix_{i+1,p}\ptl_{i+1,p})(\phi_{\vec a})\equiv c_{i,r}\phi_{\vec a}
\;\;({\rm mod}\;\sum_{s=1}^iU_s),\qquad
c_{i,r}\in\mbb{C},\eqno(3.24)$$ for $2\leq i\leq n-1$ and $1\leq
r\leq i$. Since for $2\leq i\leq n-2$,
$$E_{i,n}=[E_{i,i+1},[E_{i+1,i+2},\cdots,[E_{n-2,n-1},E_{n-1,n}]\cdots]],\eqno(3.25)$$
we have
$$d_{i,n}=[d_i,[d_{i+1},\cdots [d_{n-2},d_{n-1}]\cdots]].\eqno(3.26)$$
Thus
$$d_{i,n}(z)=0\qquad\for\;\;2\leq i\leq n-1.\eqno(3.27)$$
By (3.11), (3.16) and (3.22)-(3.27), we obtain
$$d_{i,n}(z)=[(\bar{\lmd}_i-c_{n-1,i})\ptl_{n,i}(f)+\sum_{q=1}^{i-1}x_{i,q}\ptl_{n,q}(f)]
x_{2,1}^{a_1}\phi_{\vec a}-\sum_{j=i+1}^{n-1}\ptl_{j,i}d_{j,n}(z)
\equiv 0\;\;({\rm mod}\;U)\eqno(3.28)$$ for $i=2,3,...,n-2$ and
$$d_{n-1}(z)=(\sum_{q=1}^{n-2}x_{n-1,q}\ptl_{n,q}(f))x_{2,1}^{a_1}\phi_{\vec a}\equiv 0
\;\;({\rm mod}\;U).\eqno(3.29)$$ Since the constraint on
$d_{r,n}(z)\equiv 0\;(\mbox{mod}\:U)$ for $r\geq 2$ implies
$\ptl_{r,s}d_{r,n}(z)\equiv 0\;(\mbox{mod}\:U)$ for $s=1,2,...,r-1$
by (3.23), (2.28) is equivalent to
$$d_{i,n}(z)=[(\bar{\lmd}_i-c_{n-1,i})\ptl_{n,i}(f)+\sum_{q=1}^{i-1}x_{i,q}\ptl_{n,q}(f)]
x_{2,1}^{a_1}\phi_{\vec a} \equiv 0\;\;({\rm mod}\;U)\eqno(3.30)$$
for $i=2,3,...,n-2$.

Expressions (3.29) and (3.30) give
$$\sum_{q=1}^{i-1}x_{i,q}\ptl_{n,q}(f)+(\bar{\lmd}_i-c_{n-1,i})\ptl_{n,i}(f)=0
\qquad\for\;\;i=2,...,n-2,\eqno(3.31)$$
$$\sum_{q=1}^{n-2}x_{n-1,q}\ptl_{n,q}(f)=0.\eqno(3.32)$$
We view
$$\ptl_{n,1}(f),\;\ptl_{n,2}(f),\;\cdots,\;\ptl_{n,n-2}(f)\;\;\mbox{ as unknowns}.\eqno(3.33)$$
Then the coefficient determinant of the system (3.31) and (3.32) is
\begin{eqnarray*}\hspace{2cm}& &\left|\begin{array}{cccc}x_{2,1}& \lmd_2-c_{n-1,2}& &\\
x_{3,1}& x_{3,2}&\ddots&\\ \vdots&\ddots&\ddots
&\lmd_{n-2}-c_{n-1,n-2}\\ x_{n-1,a}&\cdots &x_{n-1,n-3}&
x_{n-1,n-2}\end{array}\right|\\
&=&\prod_{p=2}^{n-1}x_{p,p-1}+g(x_{2,1},x_{3,2},...,x_{n-1,n-2})\not\equiv
0,\hspace{4.6cm}(3.34)\end{eqnarray*} where
$g(x_{2,1},x_{3,2},...,x_{n-1,n-2})$ is a polynomial of degree $n-3$
in $\{x_{2,1},x_{3,2},...,x_{n-1,n-2}\}$ over ${\cal A}_0$ (cf.
(2.32)). Therefore,
$$\ptl_{n,q}(f)=0\qquad\for\;\;q=1,2,....,n-2.\eqno(3.35)$$

Based on our calculations in (3.23)-(3.25), we can prove by
induction that
$$\ptl_{q+r,q}(f)=0\qquad\for\;\;1\leq q\leq n-2.\;2\leq r\leq n-q.\eqno(3.36)$$
So  $f$ is  a constant.

Suppose that $z$ is any solution of the system (3.4) in ${\cal
A}_1$. By (2.34) and (3.2), it can be written as
$$z=\sum_{\vec j\in\mbb{N}^{n-1}}\;\sum_{i=1}^pf_{\vec a^i-\vec j}
x_{2,1}^{a^i_1-j_1}\phi_{\vec a^i-\vec
j}\qquad\mbox{with}\;\;f_{\vec j}\in {\cal A}_0. \eqno(3.37)$$ Let
$${\cal S}=\{\vec b\in\mbb{C}^{\:n-1}\mid f_{\vec b}\neq 0;\;f_{\vec b+\vec j}=0\;
\mbox{for all}\;\vec 0\neq\vec j\in\mbb{N}^{\:n-1}\}.\eqno(3.38)$$
The above arguments show that
$$\{f_{\vec b}\mid \vec b\in {\cal S}\}\qquad \mbox{are constants}\eqno(3.39)$$
(cf. the key equations (3.29) and (3.30)). Since $\sum_{\vec b\in
{\cal S}}f_{\vec b}x_{2,1}^{b_1}\phi_{\vec b}$ is a solution of the
system (3.4), so is $z-\sum_{\vec b\in {\cal S}}f_{\vec
b}x_{2,1}^{b_1}\phi_{\vec b}$. By induction, we prove the Claim\psp

 To solve the system (2.27) in ${\cal A}_1$, we only need to
consider the solutions of the form
 $z=x_{2,1}^{a_1}\phi_{\vec a}$ with $\vec a\in\mbb{C}^{\:n-1}$ by the above claim, because
  (2.27) is a weighted system. Note
$$d_1=(\lmd_1-\sum_{j=2}^nx_{j,1}\ptl_{j,1}+\sum_{j=3}^nx_{j,2}\ptl_{j,2})\ptl_{2,1}
-\sum_{j=3}^nx_{j,2}\ptl_{j,1}.\eqno(3.40)$$ Denote
$$\hat a_{1,r}=\sum_{p=r}^{n-1}a_i-\sum_{p=2}^{n-r}(p-1)(\lmd_p+1)-(n-r)\sum_{q=n-r+1}^{n-1}
(\lmd_q+1),\;\;\hat a_{1,1}=a_1+\hat a_{1,2}\eqno(3.41)$$ for
$r=2,3,...,n-1$,
$$\hat a_{2,r}=a_r-\sum_{p=1}^{r-2}(\lmd_{n-p}+1)\qquad \for\;\;r=2,3,...,n-1.\eqno(3.42)$$
and
$$\td{a}=\sum_{i=2}^{n-1}a_i-\sum_{p=3}^{n-1}(p-2)(\lmd_p+1).\eqno(3.43)$$
Letting $x_{p,q}=0$ for $1\leq q\leq p-2\leq n-2$ in
$$d_1(z)=[(\lmd_1-\sum_{j=2}^nx_{j,1}\ptl_{j,1}+\sum_{j=3}^nx_{j,2}\ptl_{j,2})\ptl_{2,1}
-\sum_{j=3}^nx_{j,2}\ptl_{j,1}](x_{2,1}^{a_1}\phi_{\vec
a})=0,\eqno(3.44)$$ we get
$$\hat a_{1,1}(\lmd_1+1-\hat a_{1,1}+\td{a})-\sum_{r=2}^{n-1}\hat a_{2,r}\hat a_{1,r}=0\eqno(3.45)$$
by (2.36) and (3.2).

Suppose  $n>3$. We take $x_{p,\:q}=0$ for $1\leq q\leq p-2\leq n-2$
in
$$\ptl_{n,1}d_1(z)=\ptl_{n,1}[(\lmd_1-\sum_{j=2}^nx_{j,1}\ptl_{j,1}+\sum_{j=3}^nx_{j,2}\ptl_{j,2})
\ptl_{2,1}-\sum_{j=3}^nx_{j,2}\ptl_{j,1}](x_{2,1}^{a_1}\phi_{\vec
a})=0,\eqno(3.46)$$ and obtain
\begin{eqnarray*} \hspace{2cm}& &a_{n-1}\left[\prod_{i=1}^{n-2}(a_{n-1}-i-\sum_{p=1}^i
\lmd_{n-p})\right][(\hat a_{1,1}-1)(\lmd_1-\hat a_{1,1}+\td{a})\\ &
&-\sum_{r=2}^{n-2}\hat a_{2,r} (\hat a_{1,r}-1)-(\hat
a_{2,n-1}-1)(\hat
a_{1,n-1}-1)]=0.\hspace{3.4cm}(3.47)\end{eqnarray*} Note that
\begin{eqnarray*}& &[\hat a_{1,1}(\lmd_1+1-\hat a_{1,1}+\td{a})-\sum_{r=2}^{n-1}\hat a_{2,r}\hat a_{1,r}]
-[(\hat a_{1,1}-1)(\lmd_1-\hat a_{1,1}+\td{a})\\
&&-\sum_{r=2}^{n-2}\hat a_{2,r}(\hat a_{1,r}-1)-
(\hat a_{2,n-1}-1)(\hat a_{1,n-1}-1)]\\
&=&\lmd_1+\td{a}-\sum_{r=2}^{n-1}\hat a_{2,r}+1-\hat a_{1,n-1}
=\lmd_1+1-\hat a_{1,n-1}\\ &=& (n-1)+\sum_{i=1}^{n-1}\lmd_i-a_{n-1}.
\hspace{9.6cm}(3.48)\end{eqnarray*} By (3.45), (3.47) and (3.48), we
have
$$a_{n-1}\prod_{i=1}^{n-1}(a_{n-1}-i-\sum_{p=1}^i\lmd_{n-p})=0.\eqno(3.49)$$
Therefore,
$$a_{n-1}\in\{0,i+\sum_{p=1}^i\lmd_{n-p}\mid i=1,2,...,n-1\}.\eqno(3.50)$$

Assume $n=3$. Then
$$d_1=(\lmd_1-x_{2,1}\ptl_{2,1}-x_{3,1}\ptl_{3,1}+x_{3,2}\ptl_{3,2})\ptl_{2,1}-x_{3,2}\ptl_{3,1},
\eqno(3.51)$$
$$z=x_{2,1}^{a_1}\phi_{\vec a}=x_{2,1}^{a_1}(x_{3,2}+x_{3,1}\ptl_{2,1})^{a_2}
(x_{2,1}^{a_2-\lmd_2-1}),\eqno(3.52)$$ and (3.45) becomes
$$(a_1+a_2-\lmd_2-1)(\lmd_1+\lmd_2+2-a_1)-a_2(a_2-\lmd_2-1)=0\eqno(3.53)$$
Letting $x_{3,1}=0$ in
\begin{eqnarray*}\hspace{1.5cm}&&\ptl_{3,1}d_1(z)=\ptl_{3,1}[(\lmd_1-x_{2,1}\ptl_{2,1}-x_{3,1}\ptl_{3,1}
+x_{3,2}\ptl_{3,2})\ptl_{2,1}\\&
&-x_{3,2}\ptl_{3,1}][x_{2,}^{a_1}(x_{3,2}
+x_{3,1}\ptl_{2,1})^{a_2}(x_{2,1}^{a_2-\lmd_2-1})]=0,\hspace{4.4cm}(3.54)\end{eqnarray*}
we get
\begin{eqnarray*}\hspace{1cm}& &a_2(a_2-\lmd_2-1)(a_1+a_2-\lmd_2-2)(\lmd_1+\lmd_2+1
-a_1)\\ &
&-a_2(a_2-1)(a_2-\lmd_2-1)(a_2-\lmd_2-2)=0,\hspace{5.4cm}(3.55)\end{eqnarray*}
equivalently
$$a_2(a_2-\lmd_2-1)[(a_1+a_2-\lmd_2-2)(\lmd_1+\lmd_2+1-a_1)-(a_2-1)(a_2-\lmd_2-2)]
=0.\eqno(3.56)$$ By (3.53), we have
\begin{eqnarray*}\hspace{1cm}& &(a_1+a_2-\lmd_2-2)(\lmd_1+\lmd_2+1-a_1)-(a_2-1)
(a_2-\lmd_2-2)\\ &=& -(\lmd_1+a_2)+a_2(a_2-\lmd_2-1)
-(a_2-1)(a_2-\lmd_2-2)\\ &=&-(\lmd_1+a_2)+2a_2-\lmd_2-2\\
&=&a_2-\lmd_1-\lmd_2-2. \hspace{9.7cm}(3.57)\end{eqnarray*} Thus
(3.56) and (3.57) give
$$a_2(a_2-\lmd_2-1)(a_2-\lmd_1-\lmd_2-2)=0,\eqno(3.58)$$
which implies that (3.50) holds for any $n\geq 2$.

When $n=2$, the solution space of (2.27) is
$\mbb{C}+\mbb{C}\sgm_1(1)$ by (2.66). In general, we
 can use (3.50) to reduce the problem of solving (2.27) to $sl(n-1)$ as follows.  Denote
\begin{eqnarray*}\hspace{1.5cm}& &\Psi_i=x_{2,1}^{a_1}\eta_2^{a_2}\eta_1^{a_2-\lmd_{n-1}-1}
\cdots\eta_r^{a_r}\eta_{r-1}^{a_r-\lmd_{n-1}-1}\cdots\eta_1^{a_r-\lmd_{n-1}-\cdots-
\lmd_{n-(r-1)}-(r-1)}\\ & &\cdots
\eta_{n-2}^{a_{n-2}}\eta_{n-3}^{a_{n-2}-\lmd_{n-1}-1}
\cdots\eta_1^{a_{n-2}+3-\lmd_3-\cdots-\lmd_{n-1}-n}\\ &
&\eta_{i-2}^{-\lmd_{i-1}-1}
\eta_{i-3}^{-\lmd_{i-1}-\lmd_{i-2}-2}\cdots
\eta_1^{-\lmd_3-\cdots-\lmd_{i-1}-(i-2)}
\hspace{5.1cm}(3.59)\end{eqnarray*} for $i=1,...,n-1,n$, where we
treat
$$\eta_{i-2}^{-\lmd_{i-1}-1}\eta_{i-3}^{-\lmd_{i-1}-\lmd_{i-2}-2}\cdots
\eta_1^{-\lmd_i-\cdots-\lmd_{i-1}-(i-2)}=1\qquad\mbox{if}\;\;i=1,2.\eqno(3.60)$$
Moreover, we set
$$\psi_n=1,\;\;\psi_i=\sgm_{n-1}\sgm_{n-2}\cdots\sgm_i(1)=\eta_{n-1}^{n-i+\sum_{p=1}^{n-i}\lmd_{n-p}}\eta_{n-2}^{n-i-1+
\sum_{p=2}^{n-i}\lmd_{n-p}}\cdots\eta_i^{\lmd_i+1}(1),\eqno(3.61)$$
for $i=1,2,...,n-1$. Then $\{\psi_i\mid i=1,...,n-1,n\}$ are
solutions of (2.27) by Theorem 2.5.  Denote
$$\lmd_{n-1,n}=0,\;\;\lmd_{n-1,i}=n-i+\sum_{p=1}^{n-i}\lmd_{n-p}\qquad\for\;\;i=1,2,...,n-1.
\eqno(3.62)$$ According to (3.50),
$$a_{n-1}=\lmd_{n-1,i_{n-1}}\qquad\mbox{for some}\;\;i_{n-1}\in\{1,...,n-1,n\}.\eqno(3.63)$$
Thus,
$$z=x_{2,1}^{a_1}\phi_{\vec a}=\Psi_{i_{n-1}}\psi_{i_{n-1}}.\eqno(3.64)$$

Set
$$\lmd_j^{(n-2)}=\left\{\begin{array}{ll}\lmd_j&\mbox{if}\;\;j<i_{n-1}-1,\\ \lmd_{i_{n-1}}
+\lmd_{i_{n-1}-1}+1&\mbox{if}\;\;j=i_{n-1}-1,\\
\lmd_{j+1}&\mbox{if}\;\; i_{n-1}\leq j\leq n-2\end{array}
\right.\eqno(3.65)$$ for $j=1,2,...,n-2$. By (2.29) and (2.30),
$$h_j(\psi_{i_{n-1}})=\zeta_j (\psi_{i_{n-1}})=\lmd^{(n-2)}_j\psi_{i_{n-1}}\qquad\for \;\;j=1,2,
...,n-2.\eqno(3.66)$$ Define
\begin{eqnarray*}\hspace{1cm}& &\bar{d_i}=(\lmd^{(n-2)}_i-\sum_{j=i+1}^{n-1}x_{j,i}\ptl_{j,i}
+\sum_{j=i+2}^{n-1}x_{j,i+1}\ptl_{j,i+1})\ptl_{i+1,i}\\
&&+\sum_{j=1}^{i-1}x_{i,j}\ptl_{i+1,j}
-\sum_{j=i+2}^{n-1}x_{j,i+1}\ptl_{j,i}\hspace{7.7cm}(3.67)\end{eqnarray*}
for $i=1,2,...,n-2$ by (2.25). Then the system
$$d_i(\Psi_{i_{n-1}}\psi_{i_{n-1}})=[d_i,\Psi_{i_{n-1}}](\psi_{i_{n-1}})=0\qquad\for\;\;i=1,2,...,n-2\eqno(3.68)$$
is equivalent to the system
$$\bar{d}_i(\Psi_{i_{n-1}}(1))=0\qquad \for\;\;i=1,2,...,n-2\eqno(3.68)$$
by Lemma 2.2, (3.64) and (3.66). The system (3.68) is a version
(2.27) for $sl(n-1)$. By induction on $n$, there exist a $\sgm'\in
S_{n-1}$ and a constant $c\in\mbb{C}$ such that
$\Psi_{i_{n-1}}\psi_{i_{n-1}}=c\sgm'(\psi_{i_{n-1}})$. Thus
$z=c\sgm'(1)$ or
$z=c\sgm'\sgm_{n-1}\sgm_{n-1}\cdots\sgm_{i_{n-1}}(1)$ for some
$i_{n-1}\leq n-1$. The last conclusion follows from (2.58) and
induction.$ \qquad\Box$\psp

If
$$n-1+\sum_{p=1}^{n-1}\lmd_p\in\mbb{N}+1,\eqno(3.69)$$
then
$$\phi=\sgm_1\cdots\sgm_{n-2}\sgm_{n-1}\sgm_{n-2}
\cdots \sgm_1(1)\eqno(3.70)$$ is a polynomial and so
$\tau^{-1}(\phi)$ (cf. (2.22)) is a nontrivial singular vector in
the Verma module $M_\lmd$ (cf. (2.14)), which was obtained by
Malikov, Feigin and Fuchs [MFF]. In general, for $1\leq i<j\leq
n-1$,
$$\phi_{i,j}=\sgm_i\cdots\sgm_{j-1}\sgm_j\sgm_{j-1}\cdots\sgm_i(1)
\;\;\mbox{is a polynomial
if}\;\;j-i+1+\sum_{r=i}^j\lmd_r\in\mbb{N}.\eqno(3.71)$$ By (3.50)
and induction on $n$, we have:\psp

{\bf Corollary 3.2}. {\it The Verma module} $M_\lmd$ {\it is
irreducible if and only if}
$$j+\sum_{p=0}^{j-1}\lmd_{i+p}\not\in\mbb{N}\qquad\mbox{\it for}\;\;1\leq i\leq n-1,\;
0\leq j\leq n-i.\eqno(3.72)$$\pse

\vspace{0.6cm}

\noindent{\Large \bf References}

\hspace{0.5cm}

\begin{description}

\item[{[BGG]}] I. N. Bernstein, I. M. Gel'fand and S. I. Gel'fand,
Structure of representations generated by vectors of highest weight,
(Russian) {\it Funktsional. Anal. i Prilozhen.} {\bf 8} (1971), no.
1, 1-9.

\item[{[DGK]}] V. V. Deodhar, O. Gabber and V. G. Kac, Structure of some categories of
 representations of infinite-dimensional Lie algebras, {\it Adv. in Math.} {\bf 45} (1982),
 no. 1, 92-116.

\item[{[DL]}] V. V. Deodhar and J. Lepowsky, On multipicity in the Jordan-H\"{o}lder series of
Verma modules, {\it J. Algebra} {\bf 49} (1977), no. 2, 512-524.

\item[{[H]}] J. E. Humphreys, {\it Introduction to Lie Algebras and Representation Theory},
 Springer-Verlag New York Inc., 1972.

\item[{[J1]}] J. C. Jantzen, Zur charakterformel gewisser darstellungen halbeinfacher grunppen
 und Lie-algebrun, {\it Math. Z.} {\bf 140} (1974), 127-149.

\item[{[J2]}] J. C. Jantzen, Kontravariante formen auf induzierten Darstellungen habeinfacher
Lie-algebren, {\it Math. Ann.} {\bf 226} (1977), no. 1, 53-65.

\item[{[J3]}] J. C. Jantzen, Moduln mit einem h\"{o}chsten gewicht, {\it Lecture Note in Math.}
 {\bf 750}, Springer, Berlin, 1979.

\item[{[K1]}] V. G. Kac, Highest weight representations of infinite-dimensional Lie algebras,
 {\it Proc. Intern. Congr. Math. (Helsinki, 1978)}, pp. 299-304, Acad. Sci. Fennica,
 Helsinki, 1980.

\item[{[K2]}] V. G. Kac,  Some problems on infite-dimensional Lie algebras and their
representations, {\it Lie algebras and related topics (New
Brunswick, N.J., 1981)}, pp. 117-126. {\it Lecture Note in Math.}
{\bf 933}, Springe, Berlin-New York, 1982.

\item[{[K3]}] V. G. Kac, {\it Infinite-Dimensional Lie algebras}, 3rd Edition, Cambridge
University Press, 1990.

\item[{[KK]}] V. G. Kac and D. A. Kazhdan, Structure of representations with highest weight of infinite-dimensional Lie algebras, {\it Adv. in Math.} {\bf 34}
(1979), no. 1, 97-108.

\item[{[KP]}] V. G. Kac and D. H. Peterson, Spin and wedge
representations of infinite-dimensional Lie algebras and groups,
{\it Proc. Natl. Acad. Sci. USA} {\bf 78} (1981), 3308-3312.

\item[{[KR]}] V. G. Kac and A. Radul, Representation theory of
the vertex algebra $W_{1+\infty}$, {\it Trans. Groups} {\bf 1}
(1996), 41-70.

\item[{[L1]}] J. Lepowsky, Canonical vectors in induced modules, {\it Trans. Amer. Math. Soc.}
{\bf 208} (1975), 219-272.

\item[{[L2]}] J. Lepowsky, Existence of canonical vectors in induced modules,
{\it Ann. of Math. (2)} {\bf 102} (1975), no. 1, 17-40.

\item[{[L3]}] J. Lepowsky, On the uniqueness of canonical vectors, {\it Proc. Amer. Math. Soc.}
 {\bf 57} (1976), no. 2, 217-220.

\item[{[L4]}] J. Lepowsky, Generalized Verma modules, the Cartan-Helgason theorem, and the Harish-Chandra homomorphism, {\it J. Algebra} {\bf 49} (1977), no. 2, 470-495.

\item[{[MFF]}] F. G. Malikov,  B. L. Feigin and D. B. Fuchs, Singular vectors in Verma modules
over Kac-Moody algebras, (Russian) {\it  Funktsional. Anal. i
Prilozhen.} {\bf 20} (1986), no. 2, 25-37.

\item[{[RW1]}]  A. Rocha-Caridi and  N. R. Walllach, Projective modules over graded Lie
algebras, I, {\it Math. Z.} {\bf 108} (1982), no. 2, 151-177.

\item[{[RW2]}]  A. Rocha-Caridi and  N. R. Walllach, Highest weight modules over graded Lie
 algebras, resolutions, filtrations and character formulas,  {\it Trans. Amer. Math. Soc.}
  {\bf 277} (1983), no. 1, 133-162.

\item[{[S]}] N. N. Sapovolov, A certain bilinear form on the universal enveloping algebra of
 a complex semisimple Lie algebras, (Russian) {\it  Funktsional. Anal. i Prilozhen.} {\bf 4}
 (1972), no. 4, 65-70.

\item[{[V1]}] D,-N. Verma, Structure of certain induced representations of complex semisimple
 Lie algebras, {\it thesis, Yale University,} 1966.

\item[{[V2]}] D,-N. Verma, Structure of certain induced representations of complex semisimple
 Lie algebras, {\it Bull. Amer. Math. Soc.} {\bf 74}(1968), 160-166.

\item[{[X]}] X. Xu, Differential invariants of classical groups, {\it Duke Math. J.} {\bf 94}
(1998), 543-572.

\end{description}

\end{document}